\documentclass[24.88pt]{article}
\usepackage{latexsym,amsfonts,euscript,amssymb}

\title{Finite  groups with conjugacy classes number one greater than its same order classes
number
\thanks{Project supported by the National Natural Science Foundation(Grant No.10171074),
 and the
Natural Science Foundation of Chongqing Education
 Commission (KG051107)
 .}
\author{ Xianglin Du \\  \small School of Mathematics,  Sichuan University,
Chengdu 610064, P.R. China and \\\small Dept. of Math.
, Chongqing Three Gorge College, Chongqing 404000, P.R. China\\
\small E-mail: duxianglin2003@yahoo.com.cn\ \\and \\ Wujie Shi \\
 \small School of Mathematics, Suzhou University, Suzhou 215006, P.R. China \\
  \small E-mail: wjshi@suda.edu.cn\\
  } }
\small \date {}

\begin{document}
\maketitle

\begin{abstract} Let $k(G)$ be the number of conjugacy classes of finite groups
$G$ and $\pi_e(G)$ be the set of the orders of elements in $G$.
Then there exists a non-negative integer $k$ such that
$k(G)=|\pi_e(G)|+k$. We call such groups to be $co(k)$ groups.
This paper classifies all finite  $co(1)$ groups. They are
isomorphic to one of the following groups:
  $A_5,\,L_2(7),\,S_5,\,Z_3,\,Z_4,\,S_4,\,A_4$, $D_{10},\,Hol(Z_5)$, or $Z_3\rtimes Z_4$.

\end{abstract}

\textbf{Keywords}\,\, finite group , conjugacy class, same order
class.

\textbf{2000 MR Subject Classification}\,\, 20D60, 20D05, 20D10

\section{Introduction}

In 1988 W. Feit and G.M. Seitz[4], and J.P. Zhang[11] proved the
well-known Syskin problem[9] independently. Syskin problem is that:
if any elements of a finite group $G$ with same order are conjugate,
then $G\cong 1, Z_2$, or $S_3$. Since then several authors have
investigated the similar problems (see [7],[8] and [12]). Motivated
by Syskin problem, we consider the following problem. We define the
following equivalence relation in a group $G$: two elements $x,y$ of
$G$ are equivalent if and only if $o(x)=o(y)$. We call this
equivalent class to be the same order class of $G$. Then there
exists a non-negative integer $k$ such that the number $k(G)$ of
conjugacy classes of $G$ satisfies $k(G)=|\pi _e(G)|+k$, where
$\pi_e(G)$ is the set of the element orders of $G$.  We call such
groups to be $co(k)$ groups. What can we say about the structure of
$co(k)$ groups when $k$ is a given number? Clearly $co(0)$ groups
are the groups discussed in Syskin problem. If $k=1$, then all
$co(1)$ groups will be classified here.

The main result of this paper is: Let $G$ be a finite group. Then
$G$ is a $co(1)$ group if and only if   $G\cong A_5, L_2(7), S_5$,
  $\,S_4,\,A_4,\,D_{10},\,Z_3,\,Z_4,\,Hol(Z_5)$ or $Z_3\rtimes Z_4 $.

  For convenience,  we give some notation. We denote the order of $x$ by $o(x)$,
  $k(G)$ the number of conjugacy classes of $G$, $\pi_e(G)$ the set of the orders
  of elements of $G$, Irr(G) the set of irreducible characters of $G$. $x\sim y$ means  that $x$ is conjugate to $y$ in
  $G$ and $x\nsim y$ means that $x$ is not conjugate to $y$ in $G$.
    $x^G$ denotes the conjugacy class containing $x$. $\bar a$ denotes the element $aN$ of
  quotient group $G/N$. The symbol 1 is used for both the identity element and the identity element
  subgroup of
  $G$. $G^*$ denotes $G-1$. Lastly we call the set $\{a^{\sigma}\,|\,a\in G,\,\sigma\in
  Aut(G)\}$ to be a fusion class of $a$ in $G$.

   By the definition,
   $co(1)$ groups are the groups whose
  conjugacy classes number is one greater than its same
  order  classes number, i.e., $co(1)$ groups have one
  and only one same order class containing two conjugacy classes of
  $G$.

    \section{Preliminaries}

      The following Lemma 2.2 and its corollary are fundamental for this paper.

    \vspace{0.2cm}
    \textbf{Lemma 2.1}\quad Suppose that $G$ is a finite $co(0)$
    group. Then $G\cong 1,\,Z_2$ or $S_3$.

    Proof. See [11; Theorem 1].

 \vspace{0.2cm}
      \textbf{ Lemma 2.2}\quad Suppose that $G$ is a finite $co(k)$ group,
      and $N$ is a normal elementary Abelian $p$-subgroup of $G$. Then $G/N$
      is a  $co(i)$ group, where $0\leq i\leq k$.

        \vspace{0.2cm}
    Proof. Let $\overline{G}=G/N$. We fix $n\in \pi_e(\overline{G}),\, (p,n)=1$. Then
    there exists a maximal order $np^i$ ($i\geq 0$) among the orders of elements in $\overline{G}$.
    We denote it by $np^{m(n)},\,m(n)\geq 0$. Clearly it dependents on $n$.
     We know that the inverse images of the elements of order $np^i$ of $\overline{G}$  consist of the
      elements of order $np^i$ or $np^{i+1}$ of $G$, and the images of the elements of order $np^{i+1}$ of $G$
      consist of elements of order $np^i$ or $np^{i+1}$ of $\overline{G}$ under the natural homomorphism
     of $G\longrightarrow G/N$.

     Let $T_i(n)$ be the set of the representatives of conjugacy classes of elements of order
     $np^i$ in $G$ and $\bar T_i(n)$ be the set of the representatives of conjugacy classes
      of element order $np^i$ in $\overline{G}$.
      Let
     $$T_{i1}(n)=\{a\in T_i(n)\,|\,o(\bar a)=np^i\},\qquad T_{i2}(n)=\{a\in T_i(n)\,|
      o(\bar a)=np^{i-1}\},$$
      $$\bar T_{i1}(n)=\{\bar a\in \bar T_i(n)\,|\,o(a)=np^i\}
      ,\qquad \bar T_{i2}(n)=\{\bar a \in \bar T_i(n)\,|\,o(a)=np^{i+1}\},$$
      $$l_i(n)=|T_i(n)|,\qquad l_{i1}(n)=| T_{i1}(n)|,\qquad
       l_{i2}(n)=| T_{i2}(n)|,$$
       $$\bar l_i(n)=|\bar T_i(n)|,\qquad \bar l_{i1}(n)=|\bar T_{i1}(n)|,\qquad
         \bar l_{i2}(n)=|\bar T_{i2}(n)|.$$
         \\It is easy to see that
      $T_i(n)=T_{i1}(n)\cup T_{i2}(n),\,T_{i1}(n)\cap
      T_{i2}(n)=\phi$ and $\bar T_i(n)=\bar T_{i1}(n)\cup \bar T_{i2}(n)
      $. So we have\quad $l_i(n)= l_{i1}(n)+l_{i2}(n)$ and $\bar l_i(n)\leq \bar l_{i1}(n)+\bar
      l_{i2}(n)$. Clearly $T_0(1)=T_{01}(1)=\{1\},\,T_{02}(1)=\phi$ and
     $\bar T_{0}(1)=\bar T_{01}(1)=\bar T_{02}(1)=\{\bar 1\}$. Hence $l_0(1)=l_{01}(1)=1,\,l_{02}(1)=0$
     and $\bar l_0(1)=\bar l_{01}(1)=\bar l_{02}(1)=1$.
     \\Now we assert that

      (1) \quad $\bar l_{i1}(n) \leq l_{i1}(n),\,\, \bar l_{i2}(n)\leq
      l_{i+1,2}(n)$, and

      (2) \quad If $G$ has the elements of order $np^{m(n)+1}$ and  the maximal order of elements of $\overline{G}$ is
      $np^{m(n)}$, then\quad $\bar l_0(n)\leq \bar l_{01}(n)+ \bar l_{02}(n)-1$.

      (1) is obvious.

      Proof of (2). Because $G$ has the elements of order $np^{m(n)+1}$
      and  the maximal order of elements of $\overline{G}$ is
      $np^{m(n)}$ for fixed $n$, therefore $\overline{G}$ has an element $\bar b$ of order $n$ and the inverse image
      $b$ is of order $np$ by the natural homomorphism of $G\longrightarrow G/N$.
      Since $(n,p)=1$, $b$ can be expressed as the product of
      an element of order $n$ and an element of order $p$. Let
       $b=b_1b_2$, where $b_1b_2=b_2b_1,\,o(b_1)=n$, and $o(b_2)=p$. Hence $b^n=b_2^n\in N$.
       Since $N$ is an elementary Abelian $p$-group, $b_2\in N$. It follows that
       $\bar b=\bar b_1\bar b_2=\bar b_1$. Since $o(b)=np$ and $o(\bar
       b)=n$, we have $\bar b\in \bar T_{02}(n),$ and $\bar b_1\in \bar
       T_{01}(n)$.  Since $\bar T_0(n)= \bar T_{01}(n)\cup \bar
       T_{02}(n)$ and $\bar b=\bar b_1 \in \bar T_{01}(n)\cap \bar
       T_{02}(n)$, we have
        $$\bar l_0(n)=|\bar T_0(n)|\leq |\bar T_{01}(n)|+|\bar T_{01}(n)|-1
        =\bar l_{01}(n)+\bar l_{02}(n)-1.$$

       Now we define a function $\delta(x)$ such that
     $\delta(n)=0$ if $G$ does not contain any element of order
     $np^{m(n)+1}$, and $\delta(n)=1$ if $G$ has  elements of order $np^{m(n)+1}$.
     It is easy
        to see that if $G$ does not contain any element of order $np^{m(n)+1}$, then
        $\bar l_{m(n)2}(n)=0$.
        Hence we have
        $$\bar l_0(n)+\bar l_1(n)+...+\bar l_{m(n)}(n)\leq$$
        $$\bar l_{01}(n)+\bar l_{02}(n)+\bar l_{11}(n)+
        \bar l_{12}(n)+...+\bar l_{m(n)1}(n)+\delta(n)(\bar
        l_{m(n)2}(n)-1).$$
        Recall that
         $\bar l_{i1}(n)\leq l_{i1}(n)$ and $\bar l_{i2}(n)\leq
        l_{i+1,2}(n)$, and if $G$ has elements of order $np^{m(n)+1}$, then
        $l_0(n)\leq \bar l_{01}(n)+\bar l_{02}(n)-1$. Therefore
                $$(2.2.1)\quad \quad\quad \quad\quad\qquad \bar l_0(n)+\bar l_1(n)+...+\bar l_{m(n)}(n)\leq$$
        $$\bar l_{01}(n)+\bar l_{02}(n)+\bar l_{11}(n)+
        \bar l_{12}(n)+...+\bar l_{m(n)1}(n)+\delta(n)(\bar l_{m(n)2}(n)-1)\leq$$
       $$ l_{01}(n)+l_{12}(n)+l_{11}(n)+l_{22}(n)+l_{21}(n)+...+
        l_{m(n)2}(n)+l_{m(n)1}(n)+\delta(n)(l_{m(n)+1,2}(n)-1)$$
        $$\leq l_0(n)+l_1(n)+l_2(n)+...+l_{m(n)}(n)+\delta(n)(l_{m(n)+1}(n)-1).$$
        That is
     $$\sum_{i=0}^{m(n)}\bar l_i(n)\leq
     \sum_{i=0}^{m(n)}l_i(n)+\delta(n)(l_{m(n)+1}(n)-1).$$
     Let $n$ run through all elements of $\pi_e(\overline G)$ such that $(n,p)=1$
     ($n$ may be 1). Then
     $$(2.2.2)\qquad \sum_{n\in\pi_e(\overline G)}^{(p,n)=1}\sum_{i=0}^{m(n)}\bar l_i(n)
     \leq \sum_{n\in\pi_e(\overline G)}^{(p,n)=1}\sum_{i=0}^{m(n)}l_i(n)+
     \sum_{n\in\pi_e(\overline G)}^{(p,n)=1}
     \delta(n)(l_{m(n)+1}(n)-1).$$
     Since $\pi_e(\overline G)\subseteq \pi_e(G)$, $\delta(n)=1$ if  $np^{m(n)+1}\in \pi_e(G)$
      ($np^{m(n)+1}\notin \pi_e(\overline G)$). Hence if we put $d=\sum_{n\in\pi_e(\overline G)}^{(p,n)=1}\delta(n)$, then
     the left hand side of (2.2.2) is just the number of conjugacy classes
     $k(\overline G)$. The right hand side of (2.2.2) is $k(G)-d$. Since $G$ is a $co(k)$ group,
     so the right hand side of (2.2.2) is $|\pi_e(G)|+k-d$.
     Obviously,
      $$d=\sum_{n\in\pi_e(\overline G)}^{(p,n)=1}\delta(n)=|\pi_e(G)|-|\pi_e(\overline
      G)|.$$
     So we have $$k(\overline G)\leq |\pi_e(G)|+k-(|\pi_e(G)|-|\pi_e(\overline G)|)
     =|\pi_e(\overline G)|+k.$$Therefore $\overline G=G/N$ is a $co(i)$ group with $0\leq i\leq
     k$.

 \vspace{0.2cm}
     \textbf{ Corollary 2.3} \quad Suppose that $G$ is a finite $co(k)$ group
     and $N$ is a normal elementary Abelian $p$-subgroup. If there exists nonidentity elements
     $x_1,x_2,...,x_t\in N$, such that $x_1,x_2,...,x_t$ are not conjugate each other in $G$,
     then $G/N$ is a $co(i)$ group, with $0\leq i\leq k-t+1$.

 \vspace{0.2cm}
     Proof. Since $N$ is an elementary Abelian $p$-subgroup, we have $o(\bar a)=1$ for any $a\in N$.
     Since  $o(x_i)=p,\,i=1,2,...,t$,
     $x_1,x_2,...,x_t\in T_{12}(1)$. Therefore  $l_{12}(1)\geq t$.
     Since $\bar l_{02}(1)=1$, $\bar l_{02}(1)\leq l_{12}(1)-(t-1)$.
     In the inequality (2.2.1) of Lemma 2.2, we take $n=1$ and have
       $$ \bar l_0(1)+\bar l_1(1)+...+\bar l_{m(1)}(1)\leq$$
        $$\bar l_{01}(1)+\bar l_{02}(1)+\bar l_{11}(1)+
        \bar l_{12}(1)+...+\bar l_{m(1)1}(1)+\delta(1)(\bar l_{m(1)2}(1)-1)\leq$$
       $$ l_{01}(1)+(l_{12}(1)-(t-1))+l_{11}(1)+l_{22}(1)+l_{21}(1)+...+
        l_{m(1)2}(1)+l_{m(1)1}(1)+\delta(1)(l_{m(1)+1,2}(1)-1)$$
        $$\leq l_0(1)+(l_1(1)-(t-1))+l_2(1)+l_3(1)+...+l_{m(1)}(1)+\delta(1)(l_{m(1)+1}(1)-1).$$
       Therefore when  $n$ runs across the set $\pi_e(\overline G)$ satisfying $(n,p)=1$,
       we have
       $$k(\overline G)\leq |\pi_e(G)|+k-(t-1)-(|\pi_e(G)|-|\pi_e(\overline G)|)
     =\pi_e(\overline G)+k-(t-1).$$
     Hence $G/N$ is a $co(i)$ group with $0\leq i\leq k-t+1$.

\vspace{0.2cm}
      \textbf{Lemma 2.4}\quad Let $G$ be a finite  $co(1)$ group and $N$ be a normal elementary Abelian subgroup of $G$.
     Assume that $x,\,y\in G$, $o(x)=o(y)=m$, and  $x\nsim y$ in $G$.
     If $G/N$ is also a $co(1)$ group, $o(\bar a)=o(\bar b)=n>1$, and $\bar a \nsim \bar b$,
     then  $n=o(\bar x)$ or $n=o(\bar y)$.

     Proof. Without loss of generality, let
     $o(a)=min\{o(g)\,|\,g\in aN\}$ and $o(b)=min\{o(g)\,|\,g\in bN\}$.
      Since $o(\bar a)=o(\bar b)=n$,  $o(a)=n$ or $np$. Similarly $o(b)=n$ or $np$.
    Since $G$ and $G/N$ are $co(1)$ groups, and $x\nsim y$, it is impossible that
    both $x$ and $y$ lie in $N$ by Corollary 2.3.
    Therefore at least
    one of $o(\bar x),\,o(\bar y)$ is greater than one.

     (i) Let $o(a)=o(b)=n$. If $n=m$, then $a\sim x$ or $a\sim y$ since $a,\,x,$ and $y$ are all elements of order $n$,
     and $G$ is $co(1)$ group. It follows that
     $\bar a\sim \bar x$ or $\bar a\sim \bar y$. Therefore  $n=o(\bar a)=o(\bar x)$ or
     $n=o(\bar a)=o(\bar y)$. If $n\neq m$, then $a\sim b$, and $\bar a\sim \bar b$.
     It is contrary to the hypothesis.

     (ii) Let $o(a)=o(b)=np$. The same as (i) we may prove that
     $n=o(\bar a)=o(\bar x)$, or $n=o(\bar a)=o(\bar y)$.

     (iii) Let $o(a)=n$ and $o(b)=np$. If $n=m$, then we have  $n=o(\bar a)=o(\bar x)$, or $n=o(\bar a)=o(\bar y)$
     as the same as case (i).
     If $n\neq m$, then $a,\,b^p$ are both elements of order $n$. Hence $a\sim b^p$,
     and $\bar a\sim \bar b^p$ follows. Therefore $o(\bar a)=o(\bar b^p)=o(\bar b)$, it implies
     that $(n,p)=1$. Let $b=b_1b_2,\,o(b_1)=n$, and $o(b_2)=p$. Since  $b^n=b_2^n\in N$,
     $b_2\in N$. Therefore $\bar b=\bar b_1$. But $o(b_1)=n<o(b)=np$, it contradicts
     the choice of $b$.

\vspace{0.2cm}
 \textbf{Lemma 2.5} [Ito] \quad Suppose that $G$ is a finite group and $A$ is
an Abelian subgroup  $G$, then for any irreducible character
$\chi$, $\chi(1)\leq |G:A|$. If $A$ is an Abelian normal subgroup
of $G$ then $\chi (1)||G:A|$.

\vspace{0.2cm} Proof. See [5; Problem (5.4) and Theorem (6.15)].

\vspace{0.2cm}
 \textbf{Lemma 2.6} \quad Suppose that $G$ is a finite group,
$N\lhd G$, $k(G)=r$, and $k(G/N)=s$, then there exists
$t$($t=r-s$) distinct irreducible characters
$\chi_1,\chi_2,...,\chi_t$ of $G$ such that

$|G/N|(|N|-1)=\chi^2_1(1)+\chi^2_2(1)+...+\chi^2_t(1)$

\vspace{0.2cm}
 Proof. Since  $|$Irr(G)$|$=$k(G)=r$ and $|$Irr(G/N)$|$=$k(G/N)=s$, $G/N$ has $s$ distinct irreducible
characters $\theta_1,\theta_2,...,\theta_s$, they may be viewed as
distinct irreducible characters of $G$. Let $\chi_1,
\chi_2,...,\chi_t$ be the all other irreducible characters of $G$.
Then
$|G|=\theta^2_1(1)+\theta^2_2(1)+...+\theta^2_s(1)+\chi^2_1(1)+\chi^2_2(1)+...+\chi^2_t(1)$=
$|G/N|+\chi^2_1(1)+\chi^2_2(1)+...+\chi^2_t(1)$. Therefore
$|G/N|(|N|-1)=\chi^2_1(1)+\chi^2_2(1)+...+\chi^2_t(1)$.

\vspace{0.2cm}
  \textbf{Lemma 2.7}\quad Let $G$ be a finite non-Abelian simple
      group. If any same order elements are contained in at most
      two fusion classes, then $G$ is isomorphic to one of the
      following groups:

      (a) $A_5,\,A_6,\,A_7,\,A_8$;

       (b) $M_{11},\,M_{12},\,M_{22},\,M_{23},\,J_2,\,M_cL$;

     (c) $L_2(7),\,L_2(8),\,L_2(11),\,L_2(16),\,L_2(27),\,L_3(3),\,L_3(4),\,
     U_3(3),\,U_3(4),\,U_3(5),\,U_3(8),$

     $Sp_4(4),\,Sz(8)$.

\vspace{0.2cm}
 Proof. See [8; Theorem 1.1].

\vspace{0.2cm}
     \textbf{Lemma 2.8} \quad Let $G$ be a finite $co(1)$ group, $N\lhd
     G$ and $N=N_1\times N_2\times...\times N_t$, where each $N_i$
     is a non-Abelian simple group and isomorphic to each other,
     i=1,2,...,t. Then $t=1$.

\vspace{0.2cm}
     Proof. If $t>1$, put $a_1,\,b_1\in N_1$ and $o(a_1)\neq o(b_1)$. Since $G$ is a $co(1)$ group,
     then exists at least one  of the same order classes of order
     $o(a_1)$ or $o(b_1)$ formed by one conjugacy class from the definition of $co(1)$ group. Without loss
     of generality we may assume all elements of order $o(a_1)$ lie in same conjugacy class.
      Let $a_2\in N_2$, such that $o(a_1)=o(a_2)$. Since $a_1a_2=a_2a_1$, $o(a_1)=o(a_1a_2)$.
      Therefore there exists $g\in G$ such that $a_1a_2=a_1^g$.
      Since $a_1^g\in N_1^g$ and
      $N^g=N_1^g\times N_2^g\times ...\times N_t^g=N_1\times N_2\times...\times
      N_t$, $N_1^g=N_i,\,1\leq i\leq t$. It follows that
            $a_1a_2=a_1^g\in N_1^g=N_i$. From $a_1\in N_1,\,a_2\in
            N_2$ we have $a_1a_2\in (N_1\times N_2)\cap N_i=1, i\neq 1,2$, a contradiction. So $t=1$.

 \vspace{0.2cm}
    \textbf{Lemma 2.9}\quad Let $G$ be a finite group,  $N\lhd G$. If $N$ is
         a non-Abelian simple
         group, then $G/C_G(N)N\cong$ a subgroup of $Out(N)$.

         \vspace{0.2cm}
         Proof. Since $N$ is a non-Abelian simple group, it is easy to get from $N/C$ theorem.

          \vspace{0.2cm}

        In later proof, we will use the following lemma repeatedly.

\vspace{0.2cm}
              \textbf{Lemma 2.10}\quad  Let $G$ be a finite group,
              $N\lhd G,\,a\in N$, and $o(a)>1$.
              Suppose all elements of order $o(a)$ of $N$ lie in $m$ $N$-conjugacy classes,
              but in one $G$-conjugacy class. If the length of $N$-conjugacy classes
              is  the same, then
             $$|G:N|=m|C_G(a):C_N(a)|$$
             Proof. By the hypothesis, $|a^G|=m|a^N|$. Since
             $|G:C_N(a)|=|G:C_G(a)||C_G(a):C_N(a)|=|G:N||N:C_N(a)|$, $|a^G|=|G:C_G(a)|$, and
             $|a^N|=|N:C_N(a)|$, therefore $|a^G||C_G(a):C_N(a)|=|G:N||a^N|$.
             The equality $|a^G|= m|a^N|$ implies that $|G:N|=m|C_G(a):C_N(a)|$.

                   \vspace{0.2cm}
               \textbf{Lemma 2.11}\quad  Suppose that $G$ is an nonsolvable $co(1)$ group, and $G$ is
                not simple. Let $N$ be a minimal normal subgroup of
                $G$. If $N$ is an elementary Abelian $p$-subgroup,
                then $G/N$ is not isomorphic to $A_5,\,L_2(7)$, or $\,S_5$.

                \vspace{0.2cm}
                Proof.
                 By Lemma 2.2,  $G/N$ is a $co(0)$ or $co(1)$ group.
                If $N^*$ contains two distinct conjugacy classes of $G$, then $G/N$ is a $co(0)$ group by
                  Corollary 2.3. Furthermore $G/N\cong 1,\,Z_2,\,S_3$ by Lemma 2.1.
                    This contradicts to $G$ is nonsolvable.
                    Therefore $G/N$ is a $co(1)$ group and
                    $N^*$  consists of one conjugacy class of
                  $G$. Let $N=1\cup a^G, a \in G$ and
                  $|N|=p^{\alpha}$.  Now,
                  we prove
                  that $G/N$ is not isomorphic to one of
                $A_5,\,L_2(7)$, or $\,S_5$.

                 \vspace{0.2cm}
                 (I) $G/N\ncong A_5$. If $G/N\cong A_5$, then $(p^{\alpha}-1)\,|\,2^2.3.5$.
                   (1) If $p=5$, then $|G|=2^2.3.5.5^{\alpha}$.
                 Since $|a^G|=(5^{\alpha}-1)\,|\,|G|$, we have
                 $\alpha=1,\,|N|=5$ and $|a^G|=4$. Thus
                 $|C_G(N)|=|C_G(a)|=|G|/|a^G|=3.5^2$ and $N<C_G(N)\lhd G$,
                 contrary to that $G/N$ is simple. (2) If
                 $p\neq 5$, we assert that $G$ has no elements of
                 order $5p$. Otherwise, since $(5, p)=1$, no
                 matter what $G$ has one or two conjugacy classes of
                 order $5p$,  we can conclude that the
                 elements of order 5 of $G/N\cong A_5$
                 lie in one conjugacy
                 class, contrary to the structure of
                 $A_5$. Hence $5\,|\,|a^G|=p^{\alpha}-1$.
                 We get (a) $(p,\alpha)=(11,1),\,(31,1)
                 ,\,(61,1)$,
                 (b)  $(p,\alpha)=(2,4)$ and $|C_G(a)|=2^6$,
                 $\pi_e(G)=\{1,2,3,4,5\}$. In  (a), we
                 have $1<G/C_G(N)\leq Aut(N)=Z_{p-1}$.  It is contrary to $G/N\cong
                 A_5$. In  (b), from [2; Lemma 2] we have $G\cong A_6$
                 or $G\cong N\rtimes A_5$.
                 Clearly $G$ is not isomorphic to $A_6$($A_6$ is a $co(2)$ group).
                 We know that $A_5$ has two
                 conjugacy classes of elements of order 5, so does $G$. Hence the
                 elements of order 2 are all conjugate in $N$ too by Lemma 2.4. It is impossible.

         \vspace{0.2cm}
                 (II) $G/N\ncong L_2(7)$. If $G/N\cong L_2(7)$. (1) If $p=7$, then
                 $|a^G|=7^{\alpha}-1\,|\,|G|$. So
                 $\alpha=1,\,|C_G(N)|=|C_G(a)|=2^2.7^2$. Therefore
                 $N<C_G(N)\lhd G$, contrary to $G/N\cong
                 L_2(7)$. (2) If $p\neq 7$, with the
                 same way as in (I) (2), we can prove $G$ has not
                 elements of order $7p$. Therefore
                 $7\,|\,p^{\alpha}-1$, it follows that
                 (a) $(p,\alpha)=(29,1),\,(43,1),\,(57,1)$, (b) $(p,\alpha)=(13,2)$ and (c)
                 $(p,\alpha)=(2,3),\,|N|=2^3$. In (a), we have
                 $N<C_G(N)\lhd G$, it is contrary to $G/N\cong L_2(7)$.
                 In (b), $N$ is a Hall-normal subgroup of $G$,
                 by Schur-Zassenhaus Theorem, $N$ has a complement
                 $H\cong G/N\cong L_2(7)$. $H$ acts fixed-point-free on $N$,
                 by [6; Theorem 7.24], the Sylow 2- subgroups of $H$ are
                 cyclic or generalized quaternion groups. But by [1],
                the Sylow 2- subgroups of $H\cong L_2(7)$ are dihedral
                groups, it is a contradiction.
                 In (c), since $7\nmid
                 |C_G(N)|$,  $C_G(N)$ is a $\{2,3\}$ group and
                 therefore it is a solvable group. Clearly $G/C_G(N)$ is solvable too,
                 therefore $G$ is solvable, a contradiction.

                 \vspace{0.2cm}
                 (III) $G/N\ncong S_5$.  If $G/N\cong S_5$, then
                 $(p^{\alpha}-1)\,|\,2^3.3.5$. We get
                  (1) $(p,\,\alpha)=(2,1),\,(3,1)$, $(5,1),
                 \,(7,1),\,(11,1),\,(13,1),\,(31,1),\,(41,1),\,(61,1)$.
                 (2) $(p,
                 \alpha)= (2,2), (2,4), (3,2), (5,2)$ or $(11,2)$.
                 If (1) holds, when $(p,\alpha)=(2,1)$, then $|N|=2$ and  $N\leq
                 Z(G)$. Since $G/N\cong S_5$ has two conjugacy
                 classes of elements of order 2,  if $G$ has  elements
                 of order 8, then they are in one conjugacy class by
                 Lemma 2.4. Let $G$ has $n$ cyclic subgroups of
                 order 8, then $2^2n=n\phi(8)=|x^G|=|G|/|C_G(x)|$ where
                 $x\in G,\,o(x)=8$. As $2^3\,|\,|C_G(x)|$, we
                 have $2^5\,|\,|G|=2^4.3.5$,  a contradiction.
                  We know that $\pi_e(G/N)=\{1,2,3,4,5,6\}$,
                  therefore $\pi_e(G)=\{1,2,3,4,5,6,10\}$ or
                 $\{1,2,3,4,5,6,10,12\}$ and $k(G)-k(G/N)=1$ or 2. By Lemma
                 2.6, there exist $\chi_i\in$ Irr(G) ($i=1,2$)  such that
                 $|G/N|(|N|-1)=\chi^2_1(1)$, or
                 $|G/N|(|N|-1)=\chi^2_1(1)+\chi^2_2(1)$.
                   It is easy to
                  verify that the above
                 two diophantine equations have no solution since $|G/N|(|N|-1)=2^3.3.5$.
                 When $(p,\alpha)=(3,1)$, then
                 $|C_G(N)|=|C_G(a)|=2^2.3^2.5$. In this case, if $C_G(N)$ is
                 solvable, then $G$ is solvable, if $C_G(N)$ is
                 nonsolvable, then $C_G(N)/N\cong A_5$.   Therefore
                  $\pi_e(G)=\{1,2,3,4,5,6,15\}$ or
                  $\{1,2,3,4,5,6,9,15\}$. In the same way as in case  $(p,\alpha)=(2,1)$,  we get a contradiction. When
                  $\alpha=1,\,p=5,7,11,13,31,41,61$, we have $1<G/C_G(N)\leq
                  Aut(N)=Z_{p-1}$. So both $C_G(N)$ and $G/C_G(N)$ are
                  solvable, it follows that $G$ is solvable too, a contradiction.
              Now we consider  (2).
                  If $(p, \alpha)= (2,2)$, or $(3,2)$, then we can derive a contradiction
                  as  (1) considering
                 $Aut(N)$. If $(p, \alpha)= (2,4)$, then $|a^G|=3.5$
                 and $C_G(a)= 2^7$. It implies that $C_G(N)=P_2$, where $P_2\in
                 Syl_2(G)$. So $G/C_G(N)$ and $C_G(N)$ are both
                 solvable, it implies that $G$ is solvable, a
                 contradiction.
                 If $(p, \alpha)= (11,2)$ or $(5,2)$, we can get a
                 contradiction similar to the subcase $(p,
                 \alpha)=
                 (13,2)$ in  (II) (2) (b). The proof of Lemma 2.11 is completed.

          \section{The proof of Theorems}

Now we begin to prove the main theorems.

 \textbf{Theorem 3.1}\quad
$G$ is a solvable $co(1)$ group if and only if $G$ is isomorphic
to one of the following groups:

$S_4,\,A_4,\,D_{10},\,Hol(Z_5),\,Z_3\rtimes Z_4,\,Z_3,\,Z_4$.

Proof. It is easy to check that the groups listed in Theorem 3.1
are $co(1)$ groups.   Conversely, let $G$ be a solvable $co(1)$
group. Suppose that $G$ is a minimal counterexample. Since $G$ is
solvable,  the minimal normal subgroup $N$ of $G$ is an elementary
Abelian $p$-subgroup. By Lemma 2.2, $G/N$ is a $co(0)$ or $co(1)$
group. $N^*$  consists of at most two $G$-conjugacy classes
 since $G$ is a $co(1)$ group. That is $N=1\cup a^G$ or $N=1\cup
a^G\cup b^G$.

  We suppose firstly that  $G/N$ is a $co(0)$ group, then by Lemma 2.1,  $G/N\cong
1,\,Z_2,\,S_3$.

  (1) If $G/N\cong 1$, then $G=N$ is an Abelian simple $p$-group.
  Thus it is only possible for $G\cong Z_3$, contrary to
  that $G$ is a minimal counterexample.

  (2) If $G/N\cong Z_2$.

    (i) Let $N=1\cup a^G$. If $|a^G|=1$, then $N\leq Z(G)$, it
    follows that $G$ is Abelian. $G$ is a $co(1)$ group implies
    that $k(G)=|G|=|\pi_e(G)|+1$. The only possible is
    $G\cong Z_4$, it is contrary to that $G$ is a minimal
    counterexample. Hence $|a^G|>1$. Since
    $|a^G|=(p^{\alpha}-1)\,|\,2$,
       $p=3,\,\alpha=1$, it implies that $G=S_3$
  or $Z_6$. But both $S_3$   and $Z_6$ are not  $co(1)$  groups.

        (ii) Let $N=1\cup a^G\cup b^G$, then $|a^G|>1,\,|b^G|> 1$ since $N$ is a minimal
  normal subgroup. Therefore  $|a^G|=|b^G|=2$.
  It follows that $N\cong Z_5$ and $G\cong D_{10}$, which is contrary to that $G$ is a minimal counterexample.

  (3) If $G/N\cong S_3$.

  (i) Let $N=1\cup a^G$. Since
  $|a^G|=p^{\alpha}-1\,|\,|S_3|$, $p^\alpha-1=1,\,2,\,3,\,6$.
  Therefore $p=2,\,3,\,7$.

  (a) If $p=2,\,\alpha=1$, then $|G|=12$.
  Clearly, $G$ is non-Abelian, so $G\cong A_4,\,Z_3\rtimes Z_4$ or $S_3\times
  Z_2$. $S_3\times Z_2$ is not a $co(0)$ group, hence $G\cong A_4$ or $Z_3\rtimes
  Z_4$,  contrary to  $G$ is a minimal counterexample.

   (b) If $p=2,\,\alpha=2$, then $\pi_e(G)=\{1,2,3\}$ or $\{1,2,3,4\}$. By [2; Lemma 2, 3],
  $G\cong (Z_2\times Z_2)\rtimes S_3\cong S_4$, it is contrary to that $G$ is a minimal counterexample.

      (c) If $p=3$, then $\alpha=1$. Let $P\in Syl_3(G)$,  then $|P|=9$ and
  $P\lhd G$. Clearly, $P^*$ contains at least three conjugacy
  classes of  same element order of $G$. Hence $G$ is not a $co(1)$ group.

  (d) If $p=7$, then $\alpha=1$. By $N/C$ Theorem, $S_3\cong G/N=G/C_G(N)\leq
Aut(N)=Z_6$,
  a contradiction.

  (ii) Let $N=1\cup a^G\cup b^G$. Since $N$ is a minimal normal subgroup,
     $|a^G|>1,\,|b^G|>1$, therefore $|a^G|,\,|b^G|=2,\,3,\,6$. $|N|=p^{\alpha}$
   implies that $|a^G|=|b^G|=2,\,3,\,6$. If $|a^G|=|b^G|=2$, then
   $N\cong Z_5$. Thus $G$ has a cyclic normal subgroup of order
   15 and $G\cong Z_{15}\rtimes Z_2$. It is easy to see that there are
   $\phi(15)=8$ elements of order 15 in $G$, therefore the
   elements of order 15 lies in at least 4 conjugacy classes of
   $G$, which shows that $G$ is not a $co(1)$ group. If
   $|a^G|=|b^G|=3$, then $N\cong Z_7$ and $\pi_e(G)=\{1,2,3,7,14\}$.
   Hence  $k(G)=6,\,k(G/N)=3$. By Lemma 2.6, there exists $\chi_i\in$
   Irr(G)( By
   Lemma 2.5, $\chi_i(1)\,|\,6,\,i=1,2,3$,) such that
   $36=|G/N|(|N|-1)=\chi^2_1(1)+\chi^2_2(1)+\chi^2_3(1)$, a
   contradiction. If $|a^G|=|b^G|=6$, then $N\cong Z_{13}$ and $C_G(N)=N$. By
   $N/C$ Theorem, $S_3\cong G/N\leq Aut(N)=Z_{13}$, a contradiction.

   Now suppose  that $G/N$ is a $co(1)$ group. Since $G$ is a  minimal
  counterexample,  $G/N$ is isomorphic to one of the groups listed
  in Theorem 3.1. By Corollary 2.3,  $N=1\cup a^G$. We consider
  various cases about $p$.

  Case 1\quad $p=2$. Let $P_2\in Syl_2(G)$.

  (i) If $N\leq Z(G)$, then $N\cong Z_2$. When $G/N\cong S_4$, then $\pi_e(G)=\{1,2,3,4,6\}$
   or $\{1,2,3,4,6,8\}$. So $k(G/N)=5$ and $k(G)=6$ or $7$. By Lemma 2.6, there
   exists $\chi_i\in $Irr(G) such that
   $24=|G/N|(|N|-1)=\chi^2_1(1)$(or $\chi^2_1(1)+\chi^2_2(1)$),
   a contradiction. In the same way, we can prove that
   $G/N\ncong A_4$ and $G\ncong Z_3\rtimes Z_4$. When $G/N\cong
   Hol(Z_5)$, if $8\in \pi_e(G)$, then all elements of order 8 are in at least two
   conjugacy classes with same length. Suppose that $G$ has $n$ cyclic subgroups of order 8.
   Let $o(x)=8$, then $n\phi(8)=4n=2|x^G|=2|G|/|C_G(x)|=10$, a
   contradiction.  Therefore $\pi_e(G)=\{1,2,4,5,10\}$ and
   $k(G)=6,\,k(G/N)=5$. By Lemma 2.6, there exists $\chi\in$$Irr(G)$
   such that $20=\chi^2(1)$, a contradiction. When $G/N\cong
   D_{10}$. By Lemma 2.4, the elements of order 5 lie in two
   conjugacy classes. $N\leq Z(G)$ implies that the elements of
   order 10 are in two conjugacy classes too, contrary to the fact that
   $G$ is a $co(1)$ group. $G/N\cong Z_3$ implies that $G\cong
   Z_6$ is not a $co(1)$ group. $G/N\cong Z_4$ implies that $G$ is
   an Abelian group of order 8, not a $co(1)$ group.

 (ii) If $N\nleq Z(G)$.

  (a) $G/N\cong S_4$ or $Z_3\rtimes
  Z_4$. $|a^G|=|G|/|C_G(a)|$ implies that $|N|=4$
  and $3\nmid |C_G(a)|$. So  if $G/N\cong S_4$, then $\pi_e(G)=\{1,2,3,4\}$ or
  $\pi_e(G)=\{1,2,3,4,8\}$. The former implies  $k(G)=k(G/N)=5$, it is impossible since $G$ and $G/N$ are both
  $co(1)$ groups.
  The latter implies $k(G)=6,\,k(G/N)=5$. By Lemma 2.6, there exists a $\chi\in$
    Irr(G) such that $72=|G/N|(|N|-1)=\chi^2(1)$, a contradiction.
     If $G/N\cong Z_3\rtimes
  Z_4$,  since  $G/N$ has elements of order 6,  $G$ has  elements of order 6 or
  12. Let $x$ be an element of order 6 or 12 of $G$. Let
  $x=a_2a_3=a_3a_2$. If $o(x)=12$, then $o(a_2)=4,\,o(a_3)=3$.
  Since $G/N$ has no elements of order 12,  $a^2_2=a^g\in N$,
  it follows that $3\,|\,|C_G(a)|$, a contradiction. If $o(x)=6$,
  then $o(a_2)=2,\,a_2\notin N$, it implies that the elements of order  2
  of $G$ are in two distinct conjugacy classes. Since the elements of order 4 of $G/N$
   are in two distinct conjugacy classes, contrary to Lemma 2.4.

  (b) If $G/N\cong A_4$, then $P_2\lhd G$. $|a^G|=(2^{\alpha}-1)\,|\,|A_4|$ implies that
  $|N|=4$ and $N\leq Z(P_2)$. Since $G$ and $G/N\cong A_4$ are both $co(1)$ groups,
  $\pi_e(G)=\{1,2,3,4\}$. We know that  $G/N\cong A_4$ has  two
  conjugacy classes of elements of order 3, so has $G$ by Lemma 2.4. Let $b,\,c$ be  elements of order 3, 4 of $G$, respectively.
  By class equation, we have $48=|G|=1+|a^G|+2|b^G|+|c^G|\leq
  1+3+2.16+6=42$, a contradiction.

  (c) If $G/N\cong Hol(Z_5),\,D_{10}$ or $Z_4$. It is impossible since $|Hol(Z_5)|=20,\,|D_{10}|=10,\,|Z_4|=4$,
  and  $|a^G|=(2^{\alpha}-1)\,|\,|G/N|$.

  (d) If $G/N\cong Z_3$, it is easy to conclude that $G\cong A_4$,
  contrary to  that $G$ is a minimal counterexample.

  Case 2 \quad $p=3$.

  (a) Let $G/N\cong S_4$. Since $|N|-1=|a^G|=(3^{\alpha}-1)\,|\,|S_4|$,
   $|N|=3,\,9$. We know that  $S_4$ has two conjugacy classes of elements of order 2, hence  the elements of
  order 3 of $G$ lies in one conjugacy class by Lemma
   2.4. It follows that if $|N|=3$, then
  $P_3=\langle x\rangle,\,o(x)=9$ where $P_3\in Syl_3(G)$. Suppose that $G$ has $n$ cyclic
  subgroups of order 9. Then $n\phi(9)=6n=|x^G|=|G|/|C_G(x)|=
  3.24/9=8$, it is  a contradiction. If $|N|=9$, then a Sylow 2- subgroup $P_2$ of $G$ acts fixed-point-free
   on $N$, by [6; Theorem 7.24], $P_2$ is a cyclic or generalized quaternion
   group. Clearly $P_2$ is isomorphic to the Sylow 2- subgroups of $G/N\cong S_4$, therefore $P_2$
    is a dihedral group, a contradiction.

    (b) If $G/N\cong A_4$, then  $|N|=3,\,\pi_e(G)=\{1,2,3,6\}$ or $\{1,2,3,6,9\}$
    since $(3^{\alpha}-1)\,|\,|G|$. Therefore
    $k(G)=5$ or 6. We know that $k(G/N)=4$. So by Lemma 2.6, there
    exists $\chi_i\in$ Irr(G) such that $24=\chi^2_1(1)$ or
    $24=\chi^2_1(1)+\chi^2_2(1)$, a contradiction.

    (c) If $G/N\cong Hol(Z_5)$, then  $|N|=3$ and $15\in \pi_e(G)$ since
    $|a^G|=(3^{\alpha}-1)\,|\,|G|$.
    We know that $G/N\cong Hol(Z_5)$ has two conjugacy classes of
    elements of order 4.  Thus by Lemma 2.4,
    elements of order 15 of $G$ are in one conjugacy class.
    Suppose that $G$ has $n$ cyclic subgroups of order 15, then
    $n\phi(15)=8n=|x^G|=|G|/|C_G(x)|=4$ where $x\in G, o(x)=15$, a contradiction.

    (d) If $G/N\cong D_{10}$, we have a contradiction in the same way as in
    (c).

    (e) If $G/N\cong Z_3$,  we get a contradiction immediately since $3^{\alpha}-1=|a^G|=|G|/|C_G(a)|=3$.

    (f) If $G/N\cong Z_4$, then $G\cong Z_3\rtimes Z_4$. It is
    contrary to  that $G$ is a minimal counterexample.

    Case 3 \quad $p=5$

    (a) If $G/N\cong S_4,\,A_4,\,Z_3\rtimes Z_4$ or $D_{10}$, then $G\ncong D_{10}$ and
    $\alpha=1,2$ since  $|a^G|=5^{\alpha}-1=|G|/|C_G(a)|$.

    If
    $\alpha=1$, then $G$ has an element $x$ of order 15. Firstly, assume $G\cong S_4$ or $Z_3\rtimes Z_4$.
    Since $S_4$
    has two conjugacy classes of elements of order 2 and  $Z_3\rtimes Z_4$ has
    two conjugacy classes of elements of order 4,
     by Lemma 2.4, the elements of order 15  of $G$ are in one conjugacy class. Let $G$
    has $n$ cyclic subgroups of order 15, then
    $n\phi(15)=8n=|x^G|=|G|/|C_G(x)|=8$,  it follows that $G/N\ncong Z_3\rtimes
    Z_4$,  $G/N\cong S_4$ and $n=1$. Therefore $\langle x\rangle\lhd G$, it follows that
     $G/N\cong S_4$ has a normal
     subgroup of order 3, a contradiction. Now assume $G/N\cong A_4$, then
     $\pi_e(G)=\{1,2,3,5,15\}$, $k(G)=6,\,k(G/N)=4$. By Lemma 2.6,
     there exist $\chi_1,\,\chi_2\in $Irr(G) such that
     $48=|G/N|(|N|-1)=\chi^2_1(1)+\chi^2_2(1)$. Since $G$ has an Abelian subgroup of
     order 15, we have $\chi_i(1)\leq
     4$ by  Lemma 2.5, a contradiction.

     If  $\alpha=2$, then
     $24=|a^G|=|G|/|C_G(a)|$, therefore $G/N\cong S_4$.
     Clearly, a Sylow 2- subgroup
     $P_2$ of $G$ acts fixed-point-free on $N$, so we have a
     contradiction similar to Case 2($p=3$),(A).

      (b) If $G/N\cong Hol(Z_5)$. Since $|Hol(Z_5)|=20$ and  $
      Hol(Z_5)$ has  two conjugacy classes of elements of order 4 ,  by Lemma 2.4,
      elements of order 5 of $G$ are in one conjugacy class.
      Therefore $G$ has elements of order 25, and they are in one
      conjugacy class. Let $x\in G,\,o(x)=25$, then  $\phi(25)=20\,|\,|x^G|$.
      It is impossible since $|x^G|=|G|/|C_G(x)|\leq
      4$, .

       (c) If $G/N\cong Z_3$, $Z_4$.
       $G/N\cong Z_3$ implies that $N\cong Z_5$ and $G\cong Z_{15}$, not a $co(1)$
       group.
         $G/N\cong Z_4$ implies that $N\cong Z_5$
       and  $G\cong Hol(Z_5)$, it is contrary to that $G$ is a minimal
       counterexample.

       Case 4\quad $p\neq 2,3,5$.

       Let $G/N$ be any group listed in
       Theorem 3.1, then $(|N|,|G/N|)=1$. By Schur-Zassenhaus
       Theorem, $N$ has a complement $H\cong G/N$ such that $G=NH$.  Since
       $C_G(a)=NC_H(a)$,  $|a^G|=|G|/|C_G(a)|=|H|/|C_H(a)|$.
       We know that $N=1\cup a^G$, hence $|N|-1=p^{\alpha}-1=|a^G|=|H|/|C_H(a)|$.
       $|H|=|G/N|\leq 24$ and $p\geq 7$
       implies that $\alpha =1,\,N=\langle a\rangle$ and $|C_H(a)|\leq 4$. Hence $NC_H(a)$ is an
       Abelian subgroup. By Lemma 2.6, there exist
       $\chi_i\in$Irr(G)(i=1, 2, ... , n) such that
        $|G/N|(|N|-1)=|H|^2/|C_H(a)|=\chi^2_1(1)+...+\chi^2_n(1)$,
        where $n=k(G)-k(G/N)=|\pi_e(G)|+1-(|\pi_e(G/N)|+1)
        =|\pi_e(G)|-|\pi_e(G/N)|=|\pi_e(C_H(a))|$.
        By Lemma 2.5, $\chi_i(1)\leq
        |G:NC_H(a)|=|H|/|C_H(a)|$. Since
         $|C_H(a)|\leq 4$,  $n\leq
        3$.

         If $n=1$, then
        $|C_H(a)|=1$, it follows that $|N|-1=p-1=|H|$. Therefore   we have
         $|H|=10,\,12$ and $H\cong G/N\cong D_{10},\,A_4$ or $Z_3\rtimes  Z_4$.
         By $N/C$ theorem, $G/C_G(N)=H\leq
        Aut(N)=Z_{p-1}$, a contradiction.

        If $n=2$, then $|H|^2/|C_H(a)|=\chi^2_1(1)+\chi^2_2(1)\leq
        2|H|^2/|C_H(a)|^2$. It implies that $|C_H(a)|=2$. Since $C_H(a)\lhd
        H$,  $C_H(a)\leq Z(H)$. Therefore $H\cong Z_3\rtimes
        Z_4$, it follows that    $S_3\cong H/C_H(a)=G/C_G(N)\leq
        Aut(N)=Z_{p-1}$, a contradiction.

        If $n=3$,  then $C_H(a)=Z_4$  since $|C_H(a)|\leq 4$. Therefore
        $|H|^2/4=|H|^2/|C_H(a)|=\chi^2_1(1)+\chi^2_2(1)+\chi^2_3(1)\leq
        3|H|^2/|C_H(a)|^2=3|H|^2/16<|H|^2/4$, it is a finial contradiction.

        The proof of Theorem 3.1 is completed.

\vspace{0.2cm}

 \textbf{Theorem 3.2}\quad If $G$ is a finite nonsolvable simple
group, then $G$ is
    a $co(1)$ group if and only if  $G\cong A_5,\,L_2(7)$ or $G\cong S_5$.

             \vspace{0.2cm}
             Proof. The elements of order 5 of $A_5$ are in two distinct conjugacy classes,
             the other same order elements are conjugate. Hence $A_5$ is
             a $co(1)$ group. In the same way we can verify that $L_2(7)$
             and $S_5$ are  $co(1)$ groups.

              Conversely, if $G$ is a simple $co(1)$ group, by Lemma 2.7,
             $G$ is isomorphic to one of the
             groups listed in Lemma 2.7.  Check Atlas[3], we know
             that $G\cong A_5$ or $G\cong L_2(7)$.

             If $G$ is an nonsolvable $co(1)$ group, and $G$ is not simple,
             we will prove $G\cong S_5$.

      Suppose that the result is not true. Let $G$ be a minimal counterexample.
          Let $N$ be a minimal normal subgroup of $G$, then $N$ is an elementary Abelian group,
          or the direct products of isomorphic non-Abelian simple groups.

 \vspace{0.2cm}
  \textbf{Step 1\quad $N$ is a non-Abelian simple group and
  $C_G(N)=1$}.

  Proof. If $N$ is an elementary Abelian group, by Lemma 2.2,
  $G/N$ is a $co(0)$ or $co(1)$ group. $G$ is
  nonsolvable implies that $G/N$ is nonsolvable. If $G/N$ is a $co(0)$
  group, by Lemma 2.1, $G/N\cong 1,\,Z_2,\,S_3$,   contrary to  $G$ is
  nonsolvable. Hence $G/N$ is a nonsolvable $co(1)$ group. If $G/N$ is simple,
  by above proof,  $G/N\cong A_5$
  or $ L_2(7)$; if $G/N$ is not simple,  $G/N\cong S_5$ since $G$
  is a minimal counterexample. It contradicts to  Lemma 2.11. Therefore $N$ is the direct product
   of isomorphic non-Abelian simple groups.
  By Lemma 2.8, $N$ is a non-Abelian simple group. If $C_G(N)>1$, let $K$ be a
  minimal normal subgroup of $G$ contained in $C_G(N)$,
  then  $K$ must be non-Abelian simple. Let $a\in N,\,b\in
  K,\,o(a)=o(b)=2$, then $a,\,b,\,ab$ are all elements of order 2.
  Since $NK=N\times K$,  $a$ is not conjugate to $b$. $G$ is a
  $co(1)$ group implies that $ab$ is conjugate to $a$ or $b$, it follows
  that $N\cap K>1$, a contradiction. Therefore $C_G(N)=1$.

 \vspace{0.2cm}
  \textbf{Step 2\quad$N\cong A_5$}.

           Proof.
          By step 1,  $N$ is non-Abelian simple and $C_G(N)=1$. Hence by
          Lemma 2.9,  $G/N\leq Out(N)$.
            Since $G$ is a $co(1)$
           group and Inn(G) may be viewed as subgroup of Aut(N),  the same order
           elements of $N$ belong to at most two fusion
           classes. By Lemma 2.7,  $N$ is isomorphic to one of the groups listed in
           Lemma 2.7.

           If $Out(N)$ is a cyclic group of prime order, then
           $G/N=Out(N)$ and $G$ is an automorphic extension
           of $N$. In this case  we can check out the
           $\pi_e(G)$ and $k(G)$ directly from Atlas[3]. By this way we can prove that
           $G$ is not a $co(1)$ group when
           $N\cong A_7,\,A_8,\,M_{11},\,M_{12},\,M_{22},\,M_{23},\,J_2,\,M_cL$,
      $L_2(7),\,L_2(8),\,L_2(11),\,\,L_3(3),\,
     U_3(3),\,\,Sz(8)$.
     Therefore we only need to consider $N\cong A_5,\,A_6,\,L_2(16)$, $L_2(17),\,L_3(4),\,
     U_3(4),\,U_3(5),\,U_3(8),\,S_4(4)$. (See table 1. From
     Atlas[3])

      \vspace{2cm}
           table 1
           $$\begin{array}{|l|l|l|l|l|l|l|l|l|}
           \hline
           2&2^2&4&6&2\times S_3&4&S_3&3\times S_3&4\\\hline

           A_5&A_6&L_2(16)&L_2(27)&L_3(4)&U_3(4)&U_3(5)&U_3(8)&S_4(4)
           \\\hline
           1A&1A&1A&1A&1A&1A&1A&1A&1A\\

           2A(4)&2A(8)&2A(16)&2A(28)&2A(64)&2A(320)&2A(240)&2A(536)&2A(840)\\
           3A(3)&3A(9)&3A(15)&3A(27)&3A(9)&3A(15)&3A(36)&3A(1512)&2B(840)\\
           5A(5)&3B(9)&5A(15)&B^{**}(27)&4A(16)&4A(16)&4A(8)&B^{**}(1512)&2C(256)\\
           B^*(5)&4A(4)&B^*(15)&7A(14)&4B(16)&5A(300)&5A(250)&3C(81)&3A(180)\\
            &5A(5)&15A(15)&B^*3(14)&4C(16)&B^{**}(300)&5B(25)&4A(64)&3B(180)\\
            &B^*(5)&B^{**}(15)&C^*2(14)&5A(5)&C^*2(300)&5C(25)&4B(64)&4A(32)\\
             & &C^*2(15)&13A(13)&B^*(5)&D^*3(300)&5D(25)&4C(64)&4B(32)\\
             & &D^*8(15)&B^*3(13)&7A(7)&5E(25)&6A(12)&6A(24)&5A(300)\\
             & &17A(17)&C^*4(13)&B^{**}(7)&F^*(25)&7A(7)&B^{**}(24)(7)&B^*(300)\\
             & &B^*4(17)&D^*5(13)& &10A(20)&B^{**}(7)&7A(21)&5C(300)\\
             & &C^*2(17)&E^*2(13)& &B^{**}(20)&8A(8)&B^*2(21)&D^*(300)\\
             & &D^*8(17)&F^*6(13)& &C^*7(20)&B^{**}(8)&C^*4(21)&5E(25)\\
             & &E^*6(17)&14A(14)& &D^*3(20)&10A(10)&9A(27)&6A(12)\\
             & & F^*2(17)&B^*3(14)& &13A(13)& &B^*2(27)&6B(12)\\
             & & G^*5(17)&C^*5(14)& &B^{**}(13)& &C^*4(27)&10A(20)\\
             & & H^*3(17)& & &C^*5(13)& &19A(19)&B^*(20)\\
             & & & & &D^*8(13)& &B^{**}(19)&10C(20)\\
             & & & & &15A(15)& &C^{**}2(19)&D^*(20)\\
             & & & & &B^{**}(15)& &D^*2(19)&15A(15)\\
             & & & & &C^*2(15)& &E^*4(19)&B^*(15)\\
             & & & & &D^*8(15)& &F^{**}(19)&15C(15)\\
             & & & & & & &21A(21)&D^8(15)\\
             & & & & & & &B^{**}(21)&17A(17)\\
             & & & & & & &C^{**}2(21)&B^*2(17)\\
             & & & & & & &D^*2(21)&C^*3(17)\\
             & & & & & & &E^*4(21)&D^*6(17)\\
             & & & & & & &F^{**}4(21)& \\

             \hline
           \end{array}
           $$

          ( Note: in the table, the first row denotes the Out(N), where $N$ is a non-Abelian
           simple group listed in the table. $nA(m)$ denotes the  conjugacy class of elements of order $n$.
           $m$ denotes the order of centralizer of a representative of the conjugacy
            class.)

     Case 1\quad Let $N\cong A_6$. Since $\pi_e(N)=\{1,2,3,4,5\}$ and
     $Out(N)=Z_2\times Z_2$,  $G/N\cong Z_2$, or  $G/N\cong Z_2\times Z_2$.
     If  $G/N\cong Z_2$, then we can know from Altas[3] that $G$ is
     not a $co(1)$ group. Therefore  $G/N\cong Z_2\times Z_2$. Let
     $\bar x_i(i=1,2,3)$ be the distinct elements of order 2 of
     $G/N$, then $(\bar x_i)^2\in N$. So $o(x_i)=2,4,6,8,10$.
     Since $3,5\in \pi_e(N)$, if $o(x_i)=6,10$, we can choose the representatives
     of $\bar x_i$ such that $o(x_i)=2$. Hence we may assume $o(x_i)=2,4,8$. If
     $o(x_1)=2$, since $2,4\in \pi_e(N)$,  the elements of order 2 of $G$ lie in two
     conjugacy classes. So it is only possible that
     $o(x_2)=o(x_3)=8$. $G$ is a $co(1)$ group implies that $x_2\sim
     x_3$, it follows that $\bar x_2\sim \bar x_3$, a contradiction.
     Similarly, if $o(x_1)=4$, we also have a contradiction.
     Therefore $o(x_i)=8,\,i=1,2,3$. Clearly, $x_1,\,x_2,\,x_3$
     are not conjugate to each other, which is contrary to that $G$ is a $co(1)$
     group.

     Case 2\quad  Let $N\cong L_2(16)$,  then $\pi_e(N)=\{1,2,3,5,15,17\},\,Out(N)=Z_4$.
     From table 1, we know that the elements of order 17 of $N$
     are in 8 N-conjugacy classes, so they are in at least two
     distinct G-conjugacy classes because $|G/N|\leq 4$.
     Therefore same order $17'$-elements are conjugate. Also from
     table 1, the elements of order 3 of $N$ are in one
     N-conjugacy class, they are in one G-conjugacy class too. So by Lemma 2.10, $|G:N|=|C_G(a):C_N(a)|$
     where $a\in N,\,o(a)=3$. Hence $2\,|\,|C_G(a)|$, it follows
     that $G$ has elements of order 6. Since $N$ contains all elements of order 2, 3 of $G$,
      $6\in\pi_e(N)$, a contradiction.

     Case 3\quad  Let $G/N\cong L_2(27)$.
        Since $\pi_e(N)=\{1,2,3,7,13,14\},\,Out(N)=Z_6$,
             if $G/N\cong Z_2,\,Z_3$, then we can know from
             Atlas[3]
             that $G$ is not a $co(1)$ group. Hence $G/N\cong Z_6$.
             Let $G/N=\langle \bar x\rangle$, then $x^6\in N$. Let
             $o(x)=6n,\,n\in \pi_e(N)$, then the elements of order
             $6n$ of $G$ are in two distinct conjugacy classes.
             Let $a\in N,\,o(a)=2$, from table 1, $|a^G|=|a^N|$.
              So by Lemma 2.10,  $|G:N|=|C_G(a):C_N(a)|$, it implies that  $3\,|\,|C_G(a)|$.
              Therefore $6\in \pi_e(G)$ and $6\in \pi_e(N)$ follows, a contradiction.

             Case 4\quad  Let $N=L_3(4)$, then $\pi_e(N)=\{1,2,3,4,5,7\},\,Out(N)=Z_2\times S_3$.
         Let $a,\,b\in N,\,o(a)=3,\,o(b)=5$.

        (1) If the elements of order 3 of $G$ are in two conjugacy
        classes, then  same order  $3'$-elements are conjugate.
          From table 1,  the elements of order 3 of $N$ are in one
        N-conjugacy class,  they are also in one G-conjugacy
        class. Also from table 1, the elements of order 5 of $N$
         are in two N-conjugacy classes with same length,
        but they are in one G-conjugacy class. By Lemma 2.10,
          $|G:N|=|C_G(a):C_N(a)|=2|C_G(b):C_N(b)|$.
         Therefore $2\,|\,|C_G(a)|$. Let $x\in
         C_G(a),\,o(x)=2$, then $x\in N$, it follows that $6\in
         \pi_e(N)$, a contradiction.

         (2) If  the elements of order 3 of $G$ are in one conjugacy
        class.

        (i) If the elements of order 4 of $G$ are in two
        conjugacy classes, then same as (1), we also have
          $|G:N|=|C_G(a):C_N(a)|=2|C_G(b):C_N(b)|$. Hence $6\in
         \pi_e(N)$, a contradiction.

         (ii) If the elements of order 4 of $G$ are in one
        conjugacy class, let $c,\,x\in N,\,o(c)=4,\,o(x)=2$. From table 1, we
        know that the elements of order 4 of $N$ are in three
        N-conjugacy classes with same length, the elements of order 2 of $N$
         are in one N-conjugacy class. So by Lemma 2.10,
      $|G:N|=3|C_G(c):C_N(c)|=|C_G(x):C_N(x)|$. Therefore
      $3\,|\,|C_G(x)|$. Since $x\in N$ and $N$ contains all
      elements of order 3 of $G$,  $N$ contains elements of order 6, a
      contradiction.

       Case 5 \quad Let $N=U_3(4)$, then $\pi_e(N)=\{1,2,3,4,5,10,13,15\},\,Out(N)=Z_4$.
        If $G/N\cong Z_2$, checking Atlas[3] we can know $G$ is not a
       $co(1)$ group. So $G/N\cong Z_4$. Let $G/N=\langle \bar
       x\rangle$. then $o(x)=4n,\,n\in \pi_e(N)$. Therefore the
       elements of order 4n of $G$ are in two conjugacy classes by Lemma 2.4.
              From table 1, the elements of order 3 of $N$ are in one $N$-conjugacy
              class, so they are also in one $G$-conjugacy class.
              By Lemma 2.10, $|G:N|=|C_G(a):C_N(a)|$, so $2\,|\,|C_G(a)|$ where $a\in N,\,o(a)=3$.
              Since $N$ contains all elements of order 2 of $G$, $N$ has elements of order 6,
             a  contradiction.

           Case 6\quad  Let $N=U_3(5)$,  then $\pi_e(N)=\{1,2,3,4,5,6,7,8,10\},\,Out(N)=S_3$.
        If $G/N\cong Z_2$, checking Atlas[3] we can know that $G$ is not
        a $co(1)$ group. Hence $3\,|\,|G:N|$.  From table 1, $N$ has two $N$-conjugacy classes of elements of
        order 5 with
         distinct lengths,  so the elements of order 5 of $G$ are in at least two distinct
           $G$-conjugacy classes.
           Also from table 1, the elements of order 4 of $N$ are in one
           $N$-conjugacy class, so
           they are in one $G$-conjugacy class. By Lemma 2.10,
            $|G:N|=|C_G(a):C_N(a)|,\,a\in N,\,o(a)=4$.
           Since $3\,|\,|G:N|$,   $3\,|\,|C_G(a)|$. Therefore $G$ has elements of
           order 12. $N$ contains all elements of order 3 of $G$ which implies that $12\in \pi_e(N)$,
           a contradiction.

     Case 7\quad Let $N=U_3(8)$, then $\pi_e(N)=\{1,2,3,4,6,7,9,19,21\},\,Out(N)=Z_3\times
         S_3$.          From table 1, $N$ has two  $N$-conjugacy
         classes of elements of order 3 with  distinct
         lengths. So the elements of order 3 of $G$ are at least in two
         distinct $G$-conjugacy classes. Also from table 1,  the elements of order 19 of $N$ are in
         6 distinct $N$-conjugacy classes with same length, but they are in one
         $G$-conjugacy class,
         so by Lemma 2.10, $6\,|\,|G/N|$. Let $a\in N,\,o(a)=7$,
           from table 1, $N$ has  three $N$-conjugacy
          classes  of elements of order 7 with same  length.
         Therefore  $|G:N|=3|C_G(a):C_N(a)|$ by Lemma 2.10,   it follows that $2\,|\,|C_G(a)|$ and
           $N$ has elements of order
          14, a contradiction.

           Case 8\quad  Let $N=S_4(4)$, then $\pi_e(N)=\{1,2,3,4,5,6,10,15,17\}$.
            From table 1,  $N$ has two
              $N$-conjugacy classes of elements of order 2 with
                distinct lengths, so the elements of order 2 of $G$ are in at least
              two $G$-conjugacy classes.  Similarly, the elements of order 5 of $G$ are
              in at least two $G$- conjugacy  classes, contrary to the fact that $G$ is a $co(1)$ group.
              Therefore the only possibility is that $N\cong A_5$.

 \textbf{Step 3 the finial contradiction}

      By Step 2,  we know that, $N\cong A_5$.
       Since $Out(N)=Z_2$,  $G/N\cong Z_2,\,|G|=120$.
     Therefore $G\cong S_5,\,A_5\times Z_2$, or $SL(2,5)$.
      Clearly $A_5\times Z_2$ is not a $co(1)$ group. $SL(2,5)$ has center isomorphic
                    to $Z_2$, hence it is not a $co(1)$ group.
                    Therefore $G\cong S_5$, it contradicts to the assumption that $G$ is
                    a minimal counterexample.

                    The proof of Theorem 3.3 is completed.

                    By Theorem 3.1, 3.2, we have proved the
                    following theorem:

                   \textbf{ Theorem 3.3}\quad Suppose that $G$ is a finite
group, then $G$ is a $co(1)$ group if and only if   $G\cong A_5,
L_2(7), S_5$,
  $\,S_4,\,A_4,\,D_{10},\,Z_3,\,Z_4,\,Hol(Z_5),\,Z_3\rtimes Z_4 $.


\bigskip

\begin{thebibliography}{999}







\bibitem{[1]} R. Brauer and M. Suzuki, On finite group of even
order whose 2-Sylow group is a generalized quaternion group. Proc.
Nat. Acad. Sci., 45(1959), 1757-1759.

\bibitem{[2]} R. Brandl and W.J. Shi, Finite groups whose
element orders are consecutive integers. J. Algebra. 143(1991),
388-400.

\bibitem{[3]} J.H. Conway, S.R. Norton, R.A. Parker and R.A.
Wilson, ATLAS of Finite Groups, Oxford Press, New York, 1985.

\bibitem{[4]} W. Feit and G.M. Seitz, On finite rational groups and related topics,
Illinois J. Math., 33(1989), 103-131.

\bibitem{[5]} I.M. Isaacs, Character theory of finite groups,
Academic Press, New York, 1976.

\bibitem{[6]} H. Kurzweil, Endliche
gruppen, Berlin-Heidelherg-New York, 1977.

\bibitem{[7]} C.H. Li and C.E. Praeger, The finite simple groups with
at most two fusion classes of every order. Comm. Algebra, 24(1996),
3681-3704.

\bibitem{[8])} C.H. Li, Finite groups in which every pair of elements of
        the same order is either conjugate or inverse-conjugate,
Comm. Algebra, 22(1994), 2807-2816.

\bibitem{[9]} V. D. Mazurov, The Kourovka Notebook, Unsolved problems in group
theory, AMS Trans. Ser. 2, Vol.121. AMS Procidence, RI, 1983.

\bibitem{[10]} W.J. Shi, A class of special minimal normal
subgroups(in Chinese), J. Southwest-China Teachers College,
9:3(1984), 9-13.

\bibitem{[11]}  J.P. Zhang, About Syskin problem of finite group, Sci. in China, 18:2(1988), 124-128.

\bibitem{[12]} J.P. Zhang, On finite groups of all whose elements of the same order
 are conjugate in their automorphism groups, J. Algebra, 153(1992), 22-36.










\end{thebibliography}
\end{document}